\documentclass[12pt]{amsart}
\usepackage{amsmath}
\usepackage{amsfonts}
\usepackage{epsfig}
\usepackage{fullpage,url}
\usepackage[dvips]{hyperref}
 \usepackage{amscd}   

\DeclareFontEncoding{OT2}{}{} 

\DeclareMathOperator{\group}{\GL_r\left(\mathbb{A}_F^f\right)}
\DeclareMathOperator{\goes}{\stackrel{\emph{w}}{\rightarrow}}
\DeclareMathOperator{\e}{e_{\phi}}
\DeclareMathOperator{\E}{E_{\phi}}

\DeclareMathOperator{\C}{\mathbb{C}_{\infty}}

\DeclareMathOperator{\End}{End_{K^{\sep}}}

\DeclareMathOperator{\aut}{Gal(K^{\sep}/K)}
\DeclareMathOperator{\dirac}{\overline{\delta}}

\DeclareMathOperator{\sep}{sep}

\DeclareMathOperator{\tor}{tor}
\DeclareMathOperator{\alg}{alg}
\DeclareMathOperator{\hhat}{\widehat{h}}

\DeclareMathOperator{\Frac}{Frac}

\newcommand{\GL}{\operatorname{GL}}

\newcommand{\isom}{\simeq}

\newcommand{\tensor}{\otimes}

\newtheorem{theorem}{Theorem}[section]
\newtheorem{lemma}[theorem]{Lemma}
\newtheorem{corollary}[theorem]{Corollary}

\theoremstyle{definition}
\newtheorem{definition}[theorem]{Definition}

\newtheorem{conjecture}[theorem]{Conjecture}

\theoremstyle{remark}
\newtheorem{remark}[theorem]{Remark}

\newtheorem{remarks}[theorem]{Remarks}
\title{Equidistribution for torsion points of a Drinfeld module}
\author{Dragos Ghioca}
\address{Dragos Ghioca, Department of Mathematics \& Statistics, Hamilton Hall, Room 218, McMaster University, 1280 Main Street West, Hamilton, Ontario L8S 4K1, Canada}

\email{dghioca@math.mcmaster.ca}
\begin{document}

\begin{abstract}
We prove an equidistribution result for torsion points of Drinfeld modules of generic characteristic. We also show that similar equidistribution statements provide proofs for the Manin-Mumford and the Bogomolov conjectures for Drinfeld modules.
\end{abstract}

\maketitle

\section{Introduction}
\label{se:intro}

\footnotetext[1]{2000 AMS Subject Classification: Primary, 11G09; Secondary, 11G50}

Ullmo proved in \cite{Ull} the Bogomolov Conjecture for curves embedded in their jacobians and Zhang proved in \cite{Zha} the Bogomolov Conjecture in full generality. In their proofs they obtain an equidistribution result for points of small height on an abelian variety. Bilu proves in \cite{Bilu} a similar equidistribution statement for points of small height on a power of an algebraic torus. If we restrict our attention only to torsion points of an abelian variety or of a power of the multiplicative group, the above mentioned equidistribution results provide a proof for the Manin-Mumford Conjecture for abelian varieties and, respectively, for powers of the multiplicative group. 

Using the analogy between abelian varieties and Drinfeld modules we can ask most of the questions we have for abelian varieties, also in the context of Drinfeld modules (see \cite{TS} for a proof of the Manin-Mumford Conjecture for Drinfeld modules of generic characteristic and see \cite{IMRN} for Mordell-Lang statements for Drinfeld modules of both finite and generic characteristic). Denis formulated in \cite{Den} the Manin-Mumford and the Mordell-Lang conjectures for Drinfeld modules of generic characteristic. As mentioned above, Denis questions were answered in full in the case of the Manin-Mumford problem (see \cite{TS}) and partially in the case of the Mordell-Lang problem (see \cite{IMRN}). We will formulate in Section~\ref{se:extensions} the Bogomolov Conjecture for Drinfeld modules (see Conjecture~\ref{C:Bogo}). Similarly, we can ask if the equidistribution results of Bilu, Ullmo and Zhang are valid also in the context of Drinfeld modules. In the present paper we prove an equidistribution result for torsion points of Drinfeld modules of generic characteristic (see Theorem~\ref{T:equi}). We will prove in Section~\ref{se:extensions} that possible extensions of our equidistribution result lead to a proof of the Bogomolov Conjecture and to a new proof of the Manin-Mumford Conjecture for Drinfeld modules of generic characteristic.

\section{Statement of our main result}
\label{se:statement}

We define first the notion of a Drinfeld module.

Let $p$ be a prime and let $q$ be a power of $p$. Let $A:=\mathbb{F}_q[t]$. Let $K$ be a field extension of $\mathbb{F}_q$. We fix a 
morphism $i:A\rightarrow K$. We define the operator
$\tau$ as the power of the 
usual Frobenius with the property that for every $x\in K$, $\tau(x)=x^q$. 
Then we let $K\{\tau\}$ be the ring of polynomials in
$\tau$ with coefficients from $K$ (the addition is the usual addition, while the multiplication is given by the usual composition of functions).

A Drinfeld module is a morphism $\phi:A\rightarrow
K\{\tau\}$ for which the 
coefficient of $\tau^0$ in $\phi_a$ is $i(a)$ for
every $a\in A$, and there 
exists $a\in A$ 
such that $\phi_a\ne i(a)\tau^0$. In this case, we also say that $\phi$ is defined over $K$. For every field extension $K\subset L$, the Drinfeld module $\phi$ induces an action on $\mathbb{G}_a(L)$ by $a*x:=\phi_a(x)$, for each $a\in A$.

Following the
definition 
from \cite{Goss}, we call $\phi$ a Drinfeld module of
generic characteristic 
if $\ker(i)=\{0\}$ and we call $\phi$ a Drinfeld
module of finite 
characteristic if $\ker(i)\ne \{0\}$. In the latter
case, we say that the 
characteristic of $\phi$ is $\ker(i)$ (which is a
prime ideal of $A$). If $\ker(i)=\{0\}$, then we extend $i$ to an embedding of $\Frac(A)=\mathbb{F}_q(t)$ into $K$.

We note that our definition of a Drinfeld module is not the most general one as we insist on $A$ being the ring of polynomials in one variable over $\mathbb{F}_q$. In the general case, $A$ is the ring of functions defined on a projective non-singular curve $C$, regular away from a closed point $\infty\in C$. For our definition of a Drinfeld module, $C=\mathbb{P}^1_{\mathbb{F}_q}$ and $\infty$ is the usual point at infinity on $\mathbb{P}^1$. Before stating our result, we need to introduce several technical ingredients.

\begin{definition}
\label{D:modular transcendence degree}
Let $\phi:A\rightarrow K\{\tau\}$ be a Drinfeld module. We call the modular transcendence degree of $\phi$ the smallest integer $d\ge 0$ such that a Drinfeld module isomorphic to $\phi$ is defined over a field of transcendence degree $d$ over $\mathbb{F}_q$.
\end{definition}

For the remaining of this section, unless otherwise stated, $\phi:A\rightarrow K\{\tau\}$ is a Drinfeld module of generic characteristic.

Let $v_{\infty}$ be the valuation on $F:=\mathbb{F}_q(t)$ given by the negative of the degree of any nonzero rational function, i.e. 
$$v_{\infty}\left(\frac{f}{g}\right)=\deg(g)-\deg(f)\text{ for every nonzero }f,g\in\mathbb{F}_q(t).$$
We fix an extension of $v_{\infty}$ on $K$ (we recall that $F\subset K$, as $\phi$ is a Drinfeld module of generic characteristic) and we denote it also by $v_{\infty}$. We let $K_{\infty}$ be the completion of $K$ at $v_{\infty}$. We denote by $F_{\infty}$ the completion of $F$ inside $K_{\infty}$. We fix an algebraic closure $K_{\infty}^{\alg}$ of $K_{\infty}$ and extend $v_{\infty}$ to a valuation on $K_{\infty}^{\alg}$. Finally, we let $\mathbb{C}_{\infty}$ be the completion of $K_{\infty}^{\alg}$ at $v_{\infty}$. As shown in \cite{Goss}, $\mathbb{C}_{\infty}$ is an algebraically closed, complete valued field. We let $K^{\alg}$ and $K^{\sep}$ be the algebraic and respectively, the separable closure of $K$ inside $K_{\infty}^{\alg}$. 

We define the set $\End(\phi)$ of endomorphisms of $\phi$ as the set of all $f\in K^{\sep}\{\tau\}$ such that $f\phi_a=\phi_af$, for every $a\in A$. As shown in \cite{Goss}, there exists a finite separable extension $L$ of $K$ such that each endomorphism $f\in K^{\alg}\{\tau\}$ of $\phi$ has coefficients in $L$. Moreover, $\End(\phi)$ is a finite extension of $A$ (if we identify $a\in A$ with $\phi_a\in K\{\tau\}$).

We define the torsion submodule of $\phi$ as
$$\phi_{\tor}=\{x\in K^{\alg}\mid\text{ there exists }a\in A\setminus\{0\}\text{ such that }\phi_a(x)=0\}.$$

For each nonzero $a\in A$, we let $\phi[a]=\{x\in K^{\alg}\mid \phi_a(x)=0\}$. Because $A=\mathbb{F}_q[t]$ is a PID, for each $x\in\phi_{\tor}$ there exists a unique monic polynomial $a\in A$ such that $\phi_a(x)=0$ and for every other $a'\in A$ such that $\phi_{a'}(x)=0$, then $a|a'$. We call $a$ the \emph{order} of $x$. Note that by construction, the order of a torsion point is always a monic polynomial in $t$. Also, we will always identify the greatest common divisor in $A$ of a number of polynomials by the monic generator of the principal ideal of $A$ generated by them. Finally, we note that because $\phi$ is a Drinfeld module of generic characteristic, $\phi_{\tor}\subset K^{\sep}$.

As shown by Theorem $4.6.9$ of \cite{Goss}, there exists an $A$-lattice $\Lambda\subset\C$ associated to $\phi$ (because $\phi$ has generic characteristic). Let $\e$ be the exponential function defined in $4.2.3$ of \cite{Goss} which gives a continuous (in the $v_{\infty}$-adic topology) isomomorphism $$\e:\C/\Lambda\rightarrow\C.$$
The torsion submodule of $\phi$ in $\C$ is isomorphic naturally through $\e^{-1}$ to $\left(F\tensor_{A}\Lambda\right)/\Lambda$. We let $T$ be the closure of $\phi_{\tor}\subset\C$ in the $v_{\infty}$-adic topology of $\C$. Then the restriction of $\e$ on $\left(F_{\infty}\tensor_{A}\Lambda\right)/\Lambda$ gives an isomorphism between $\left(F_{\infty}\tensor_{A}\Lambda\right)/\Lambda$ and $T$.

Let $r$ be the rank of $\Lambda$ and fix an $A$-basis $z_1,\dots,z_r$ of $\Lambda$. Then $\left(F\tensor_A\Lambda\right)/\Lambda\isom (F/A)^r$. Also, because $F=\mathbb{F}_q(t)$, then $F_{\infty}=\mathbb{F}_q((\frac{1}{t}))$. Using Proposition $4.6.3$ of \cite{Goss}, $\left(F_{\infty}\tensor_A\Lambda\right)/\Lambda$ is isomorphic to $\left(\mathbb{F}_q((\frac{1}{t}))/\mathbb{F}_q[t]\right)^r$. Then we have the isomorphism $\E:\left(\mathbb{F}_q\left(\left(\frac{1}{t}\right)\right)/\mathbb{F}_q[t]\right)^r\rightarrow T$ given by
$$\E(\gamma_1,\dots,\gamma_r):=\e(\gamma_1z_1+\dots+\gamma_rz_r),$$
for each $\gamma_1,\dots,\gamma_r\in\mathbb{F}_q\left(\left(\frac{1}{t}\right)\right)/\mathbb{F}_q[t]$.

We construct the following  group isomorphism $$\sigma:\mathbb{F}_q\left(\left(\frac{1}{t}\right)\right)/\mathbb{F}_q[t]\rightarrow\frac{1}{t}\cdot\mathbb{F}_q[[\frac{1}{t}]]\text{, given by}$$
\begin{equation}
\label{E:sigma}
\sigma\left(\sum_{i\ge -n}\alpha_i\left(\frac{1}{t}\right)^i\right)=\sum_{i\ge 1}\alpha_i\left(\frac{1}{t}\right)^i,
\end{equation}
for every natural number $n$ and for every $\sum_{i\ge -n}\alpha_i\left(\frac{1}{t}\right)^i\in\mathbb{F}_q\left(\left(\frac{1}{t}\right)\right)$ (obviously, $\sigma$ vanishes on $\mathbb{F}_q[t]$). The group $\frac{1}{t}\cdot\mathbb{F}_q[[\frac{1}{t}]]$ is a topological group with respect to the restriction of $v_{\infty}$ on $\frac{1}{t}\cdot\mathbb{F}_q[[\frac{1}{t}]]$. Hence, the isomorphism $\sigma^{-1}$ induces a topological group structure on $\mathbb{F}_q\left(\left(\frac{1}{t}\right)\right)/\mathbb{F}_q[t]$. Therefore, $\sigma$ becomes a continuous isomorphism of topological groups.
We endow $\left(\mathbb{F}_q\left(\left(\frac{1}{t}\right)\right)/\mathbb{F}_q[t]\right)^r$ with the corresponding product topology. The isomorphism $\sigma$ extends diagonally to another continuous isomorphism, which we still call $\sigma$, between $\left(\mathbb{F}_q((\frac{1}{t}))/\mathbb{F}_q[t]\right)^r$ and $G:=\left(\frac{1}{t}\mathbb{F}_q[[\frac{1}{t}]]\right)^r$. Using that $\e$ is a continuous morphism and using the definition of our topology on $\left(\mathbb{F}_q\left(\left(\frac{1}{t}\right)\right)/\mathbb{F}_q[t]\right)^r$, we conclude $\E$ is a continuous morphism. Thus 
\begin{equation}
\label{E:continuous isomorphism}
\E\sigma^{-1}:G\rightarrow T\text{ is a continuous isomorphism.}
\end{equation} 

Let $\mu$ be the Haar measure on $G$, normalized so that its total mass is $1$. Let $\nu:=\left(\E\sigma^{-1}\right)_{*}\mu$ be the induced measure on $T$ (i.e. $\nu(V):=\mu\left(\sigma\E^{-1}(V)\right)$ for every measurable $V\subset T$). Because $\mu$ is a probability measure then $\nu$ is also a probability measure. Because $\mu$ is a Haar measure on $G$ and $\E\sigma^{-1}$ is a group ismorphism, then $\nu$ is a Haar measure on $T$.

\begin{definition}
\label{D:orbit}
For each $x\in K^{\sep}$, we denote by $O(x)$ the (finite) orbit of $x$ under $\aut$.
\end{definition}

\begin{definition}
\label{D:Dirac measure}
Given $x\in K^{\sep}$, we define a probability measure $\dirac_{x}$ on $\C$ by
$$\dirac_x=\frac{1}{\# O(x)}\sum_{y\in O(x)}\delta_y,$$
where $\# O(x)$ represents, as always, the cardinality of the set $O(x)$ and $\delta_y$ is the Dirac measure on $\C$ supported on $\{y\}$.
\end{definition}

Before we can state the equidistribution result (Theorem~\ref{T:equi}), we need to define the concept of weak convergence for a sequence of probability measures on a metric space.
\begin{definition}
\label{D:weak convergence}
A sequence $\{\lambda_k\}$ of probability measures on a metric space $S$ weakly converge to $\lambda$ if for any bounded continuous function $f:S\rightarrow\mathbb{R}$, $(f,\lambda_k)\rightarrow (f,\lambda)$ as $k\rightarrow\infty$ (where $(f,\lambda)$ represents, as always, the integral of $f$ on $S$ with respect to $\lambda$). In this case we use the notation $\lambda_k\goes\lambda$.
\end{definition}

\begin{theorem}
\label{T:equi}
Let $\phi:A\rightarrow K\{\tau\}$ be a Drinfeld module of generic characteristic and of modular transcendence degree at least $2$. Assume $\End(\phi)=A$. Let $\{x_k\}$ be a sequence of distinct torsion points in $\phi$. Then $\dirac_{x_k}\goes\nu$.
\end{theorem}

\begin{remark}
\label{R:support}
If $x\in\phi_{\tor}$, then $O(x)\subset T$ and so, the measure $\dirac_x$ is supported on $T$. Therefore, the conclusion of Theorem~\ref{T:equi} should be interpreted as follows: for each $x_k$ as in Theorem~\ref{T:equi}, $\dirac_{x_k}$ is a measure on $T$ and as $k\rightarrow\infty$, the probability measures $\dirac_{x_k}$ converge weakly to the normalized Haar measure $\nu$ on $T$.
\end{remark}

\begin{remark}
\label{R:about transcendence degree}
We will explain during the proof of Theorem~\ref{T:equi} why the hypothesis on the modular transcendence degree is needed in our proof. However, we note that the modular transcendence degree of a Drinfeld module $\phi$ of generic characteristic is at least $1$, because no Drinfeld module isomorphic to $\phi$ is defined over a finite field (in that case, the Drinfeld module would be of finite characteristic).
\end{remark}

\section{Proof of the main theorem}
\label{se:proof}

We continue in this section with the notation from Section~\ref{se:statement}.
\begin{proof}[Proof of Theorem~\ref{T:equi}.]
Let $\hat{A}:=\prod_{P}A_{(P)}$ denote the profinite completion of $A$ (where $P$ runs over all the monic irreducible polynomials of $A=\mathbb{F}_q[t]$ and $A_{(P)}$ represents the completion of $A$ at the prime ideal $(P)$). We define the finite ad\`{e}les $\mathbb{A}^f_F:=F\tensor_A\hat{A}$. We let
$$\pi:\aut\rightarrow\group$$
be the natural representation on the ad\`{e}lic Tate module of $\phi$. Let $\Gamma$ be its image. The following result is Theorem $3$ of \cite{pink}.
\begin{theorem}
\label{T:pink}
Let $\phi:A\rightarrow K\{\tau\}$ be a Drinfeld module of generic characteristic. Suppose that $\End(\phi)=A$ and that $\phi$ is not isomorphic to a Drinfeld module defined over a finite extension of $F$. Then $\Gamma$ is an open subgroup of $\group$.
\end{theorem}

\begin{remark}
\label{R:explanations}
Because we will use Theorem~\ref{T:pink} in our proof we needed to impose the two extra hypothesis on $\phi$: that it has modular transcendence degree at least $2$ and its endomorphism ring equals $A$. As remarked by Pink in a \emph{Note} after the proof of his Theorem $3$ in \cite{pink}, the statement of Theorem~\ref{T:pink} is conjectured to be true without the extra assumption on the modular transcendence degree. In that case, our proof of Theorem~\ref{T:equi} would show the equidistribution result for every Drinfeld module of generic characteristic.
\end{remark}

Because our Drinfeld module $\phi$ has modular transcendence degree at least $2$, it is not isomorphic to a Drinfeld module defined over a finite extension of $F$ (otherwise, it would have modular transcendence degree $1$; see also Remark~\ref{R:about transcendence degree}). Thus $\phi$ satisfies the hypothesis of Theorem~\ref{T:pink} and so, $\Gamma$ is an open subgroup of $\group$. Hence
\begin{equation}
\label{E:finite index}
[\group:\Gamma]<\aleph_0.
\end{equation}
Using \eqref{E:finite index} we conclude there exist finitely many irreducible monic polynomials $P_1,\dots,P_l\in A$ and there exists a natural number $m$ such that 
\begin{equation}
\label{E:big group}
\prod_{i=1}^l\left(I_r+\left(\prod_{i=1}^lP_i\right)^m\GL_r\left(A_{(P_i)}\right)\right)\cdot\prod_{P\ne P_i}\GL_r\left(A_{(P)}\right)\subset\Gamma,
\end{equation}
where $I_r\in\GL_r$ is the identity matrix. Therefore, \eqref{E:big group} shows that the orbit of every torsion point of $\phi$ of order $a$ (where $a\in A$ is a monic polynomial) coprime with $\prod_{i=1}^lP_i$ consists of all the possible torsion points of order $a$. 

In general, we represent a torsion point $x$ of order $b$ (not necessarily coprime with $P_1\dots P_l$) through the isomorphism $\sigma\E^{-1}:\phi_{\tor}\rightarrow\left(\frac{1}{t}\mathbb{F}_q[[\frac{1}{t}]]\right)^r$, as 
\begin{equation}
\label{E:Q}
\left(\frac{b_1}{b},\dots,\frac{b_r}{b}\right)
\end{equation}
where 
\begin{equation}
\label{E:Qdegree}
\text{each polynomial $b_i$ has degree less than $b$}
\end{equation}
and, moreover, the greatest common divisor 
\begin{equation}
\label{E:Qgcd}
(b_1,\dots,b_r,b)=1.
\end{equation}
Then the orbit $O(x)$ of $x$ corresponds through $\sigma\E^{-1}$ to a subset $S(x)$ of the set of all possible tuples in $G=\left(\frac{1}{t}\mathbb{F}_q[[\frac{1}{t}]]\right)^r$ of the form \eqref{E:Q} satisfying \eqref{E:Qdegree} and \eqref{E:Qgcd}. Let $H\subset\aut$ be the preimage through $\pi^{-1}$ of the subgroup from \eqref{E:big group} which is contained in $\Gamma$. Then $\aut$ is a finite union of cosets of $H$. 

For each $i\in\{1,\dots,l\}$, let $\beta_i$ the exponent of $P_i$ in the prime decomposition of $b$. Let $b':=\prod_{i=1}^lP_i^{\beta'_i}$ be the monic polynomial which is the greatest common divisor between $b$ and $P:=\left(\prod_{i=1}^lP_i\right)^m$. Obviously, $\beta'_i\le\beta_i$. For each $i\in\{1,\dots,l\}$, let 
\begin{equation}
\label{E:description}
\left(\frac{\alpha_1}{P_i^{\beta_i}},\dots,\frac{\alpha_r}{P_i^{\beta_i}}\right)
\end{equation}
be the image in $S(x)$ of an element of the $P_i$-power part of $O(x)$. Then \eqref{E:big group} shows that the image through $H$ of the $P_i$-power part of the torsion point from \eqref{E:description} corresponds in $S(x)$ to all the elements of the form
\begin{equation}
\label{E:description 2}
\left(\frac{\alpha_1+P_i^{\beta'_i}q_1}{P_i^{\beta_i}},\dots,\frac{\alpha_r+P_i^{\beta'_i}q_r}{P_i^{\beta_i}}\right)
\end{equation}
for arbitrary $q_1,\dots,q_r$.

For a monic, prime polynomial $Q\in A$, different than $P_1,\dots,P_l$, let $Q^{\beta}$ be the maximal power of $Q$ dividing $b$. Then, by \eqref{E:big group}, the $Q$-power part of the orbit $O(x)$ consists of all torsion points of $\phi$ of order $Q^{\beta}$.

Let $\{x_k\}$ be a sequence of distinct torsion points of $\phi$. Because $\phi$ has generic characteristic, then $\{x_k\}\subset K^{\sep}$. Because the elements of the sequence $\{x_k\}$ are distinct, the orders $b_k\in A=\mathbb{F}_q[t]$ of each $x_k$ have the property that $\deg(b_k)\rightarrow\infty$ (the degree $\deg(b_k)$ of each $b_k$ is the degree of $b_k$ as a polynomial in $t$).

For each $k$ we let $O_k:=\sigma\E^{-1}\left(O(x_k)\right)\subset G$. For each $k$ we define by $\dirac_k$ the associated probability measure on $G$, equally supported on $O_k$:
$$\dirac_k:=\frac{1}{\#O_k}\sum_{y\in O_k}\delta_y.$$
Because $\E\sigma^{-1}:G\rightarrow T$ is a continuous isomorphism, we conclude that it suffices to show
\begin{equation}
\label{E:new conclusion}
\dirac_k\goes\mu.
\end{equation}

Let $f$ be any continuous, real valued function on $G$. Because $G$ is a totally disconnected, compact space, $f$ is a finite $\mathbb{R}$-linear combination of characteristic functions on open subsets of $G$. Hence, in order to prove \eqref{E:new conclusion}, it suffices to prove \eqref{E:new conclusion} for characteristic functions of open subsets of $G$. Thus, for each such open subset $U$, we need to show
\begin{equation}
\label{E:simplified conclusion}
\frac{\#\left(O_k\cap U\right)}{\#O_k}\rightarrow\mu(U)\text{ as }k\rightarrow\infty. 
\end{equation}

Let $U$ be an open subset of $G$. Then $U$ is of the form $$\left(a_1\left(\frac{1}{t}\right),\dots,a_r\left(\frac{1}{t}\right)\right)+\left(\frac{1}{t^{n_1+1}}\mathbb{F}_q\left[\left[\frac{1}{t}\right]\right],\dots,\frac{1}{t^{n_r+1}}\mathbb{F}_q\left[\left[\frac{1}{t}\right]\right]\right),$$
where $a_i\in\frac{1}{t}\mathbb{F}_q[\frac{1}{t}]$ is a polynomial of degree at most $n_i$ for each $i$, and $n_1,\dots,n_r$ are positive integers. Ultimately, our goal is to show
\begin{equation}
\label{E:goal}
\frac{\#(O_k\cap U)}{\#O_k}\rightarrow\mu (U)=q^{-\sum_{i=1}^r n_i}\text{, as $k\rightarrow\infty$.}
\end{equation}

Let $O$ be one of the sets $O_k$. Let $b\in\mathbb{F}_q[t]$ be the monic polynomial which is the order of $x_k$. Then all the torsion points $y\in O(x_k)$ have the same order $b$. Therefore, the elements of $O$ are of the form
\begin{equation}
\label{E:form}
\left(\frac{b_1}{b},\dots,\frac{b_r}{b}\right),
\end{equation}
where $b_i\in\mathbb{F}_q[t]$ and for each $i$, 
\begin{equation}
\label{E:degree}
\deg(b_i)<\deg(b)
\end{equation}
and moreover, the greatest common divisor 
\begin{equation}
\label{E:gcd}
(b_1,\dots,b_r,b)=1.
\end{equation}

Theorem~\ref{T:pink} will allow us to determine the proportion of elements of the form \eqref{E:form} which are in $O$. This will allow us to compute $\#(O\cap U)$, which in turn will lead to the proof of our Theorem~\ref{T:equi}.

Let $b'$ be the greatest common divisor of $b$ and $\prod_{i=1}^lP_i^m$. We let $C$ be the collection of all tuples of polynomials of the form
\begin{equation}
\label{E:coset type}
(\alpha_1,\dots,\alpha_r)\text{ with }\deg(\alpha_i)<\deg(b')\text{ for all $i$}
\end{equation}
such that for each such tuple there exists $y\in O$ of the form
\begin{equation}
\label{E:element type}
\left(\frac{\alpha_1+b'q_1}{b},\dots,\frac{\alpha_r+b'q_r}{b}\right),
\end{equation}
where for each $i$, $q_i\in\mathbb{F}_q[t]$ and 
\begin{equation}
\label{E:grad}
\deg(\alpha_i+b'q_i)<\deg(b)
\end{equation}
and
\begin{equation}
\label{E:gcd2}
\left(\alpha_1+b'q_1,\dots,\alpha_r+b'q_r,b\right)=1.
\end{equation}
If $b'=1$, then the only tuple as in \eqref{E:coset type} is $(\alpha_1,\dots,\alpha_r)=(0,\dots,0)$.

Condition \eqref{E:gcd2} shows that
\begin{equation}
\label{E:gcd3}
(\alpha_1,\dots,\alpha_r,b')=1.
\end{equation}

Clearly, $C$ is a finite set because there are finitely many tuples of the form \eqref{E:coset type} without even asking the extra condition \eqref{E:element type}. Using \eqref{E:big group} and our analysis for the action of $\aut$ on the different $Q$-power parts of $O$, we conclude that the elements of $O$ are \emph{all} the elements of $G$ of the form \eqref{E:element type}, corresponding to some tuple $(\alpha_1,\dots,\alpha_r)\in C$, and satisfying the conditions \eqref{E:grad} and \eqref{E:gcd2}.

We fix a tuple $(\alpha_1,\dots,\alpha_r)\in C$. We count the number of elements of $O$ which are of the form \eqref{E:element type} for this tuple $(\alpha_1,\dots,\alpha_r)$. As explained above, because of \eqref{E:big group} we need to count the number of all elements $G$ of the form \eqref{E:element type}, corresponding to this fixed tuple $(\alpha_1,\dots,\alpha_r)$, and satisfying \eqref{E:grad} and \eqref{E:gcd2}. 

We define the M\"{o}bius function $\mu$ on the set of all monic polynomials in $\mathbb{F}_q[t]$ by
$$\mu(1)=1,$$
$$\mu(Q_1Q_2\dots Q_n)=(-1)^n\text{ if $Q_1,\dots,Q_n$ are distinct irreducible, non-constant polynomials,}$$
$$\mu(f)=0\text{ if $f$ is not squarefree.}$$
In this proof, the letter $\mu$ also appears as denoting the measure on $G$. This should not be confused with the above defined M\"{o}bius function, as the measure $\mu$ is always evaluated on subsets of $G$, while the M\"{o}bius function $\mu$ is always evaluated on monic polynomials.

It is immediate to see that for every nonzero polynomial $f\in\mathbb{F}_q[t]$ (monic or not)
\begin{equation}
\label{E:mobius}
\sum_{g|f}\mu(g)=1\text{ if $\deg(f)=0$ and it is $0$ otherwise}.
\end{equation}
Of course, $g$ in \eqref{E:mobius} is a monic polynomials and in general, when we will sum over divisors of a polynomial $f$, we will always include only the monic polynomials dividing $f$.

Hence, in order to count the number of elements of $O$ which are of the form \eqref{E:element type} for the fixed tuple $(\alpha_1,\dots,\alpha_r)$, we compute
\begin{equation}
\label{E:sum 1}
\sum_{\substack{q_1,\dots,q_r\\\deg(\alpha_1+b'q_1)<\deg(b),\dots,\deg(\alpha_r+b'q_r)<\deg(b)}}\left(\sum_{d|(\alpha_1+b'q_1,\dots,\alpha_r+b'q_r,b)}\mu(d)\right).
\end{equation}
Using \eqref{E:mobius}, we obtain that the inner sum in \eqref{E:sum 1} equals $1$ if and only if the greatest common divisor $(\alpha_1+b'q_1,\dots,\alpha_r+b'q_r,b)=1$, otherwise the inner sum equals $0$.
Changing the order of summation in \eqref{E:sum 1} we obtain
\begin{equation}
\label{E:sum 2}
\sum_{d|b}\mu(d)\cdot\left(\sum_{\substack{q_1,\dots,q_r\\\deg(\alpha_1+b'q_1)<\deg(b),\dots,\deg(\alpha_r+b'q_r)<\deg(b)\\d|(\alpha_1+b'q_1),\dots,d|(\alpha_r+b'q_r)}}1\right).
\end{equation}
We evaluate the inner sum taken into account that for its computation, $\alpha_1,\dots,\alpha_r$ and $b$, $b'$ and $d$ are all fixed. 
We also take into account that if $d$ and $b'$ are not coprime, then the inner sum is $0$ as shown by \eqref{E:gcd3}. On the other hand, if $d$ and $b'$ are coprime, then each congruence
\begin{equation}
\label{E:congruence}
\alpha_i+b'q_i\equiv 0\text{ (mod $d$) has one incongruent solution $q_i$ modulo $d$.}
\end{equation}
For each $i$, let $s_i$ be the unique solution to the congruence \eqref{E:congruence} with $\deg(s_i)<\deg(d)$. Then all the solutions $q_i$ to the congruence \eqref{E:congruence}, which we count in the inner sum from \eqref{E:sum 2}, are of the form 
\begin{equation}
\label{E:q'}
q_i=s_i+dq_i',
\end{equation}
for some polynomials $q_i'$ such that 
\begin{equation}
\label{E:degree equation}
\deg(\alpha_i+b'q_i)<\deg(b). 
\end{equation}
Using \eqref{E:q'} we get 
\begin{equation}
\label{E:prima parte}
\alpha_i+b'q_i=\alpha_i+b's_i+b'dq_i'.
\end{equation}
Because the inner sum of \eqref{E:sum 2} is nonzero only if $(d,b')=1$ and because $d|b$, we need to estimate the inner sum of \eqref{E:sum 2} only when $d|\frac{b}{b'}$ (and even then, the inner sum still might be $0$). So, for us, $\deg(d)\le\deg(b)-\deg(b')$. This shows
$$\deg(\alpha_i+b's_i)\le\max\{\deg(\alpha_i),\deg(b')+\deg(s_i)\}<\max\{\deg(b'),\deg(b')+\deg(d)\}\le\deg(b).$$
Therefore, for each $i$ and for each polynomial $q'_i$ of degree less than $\deg(b)-\deg(b'd)$, the degree of \eqref{E:prima parte} is less than the degree of $b$. Hence, we have $q^{(\deg(b)-\deg(b'd))r}$ choices for tuples $(q_1',\dots,q_r')$. So, we compute the inner sum in \eqref{E:sum 2} and obtain
\begin{equation}
\label{E:inner sum}
q^{r\left(\deg(b)-\deg(b'd)\right)}=q^{r\left(\deg(b)-\deg(b')\right)-r\deg(d)},
\end{equation}
if $(d,b')=1$, while if $(d,b')\ne 1$, the inner sum in \eqref{E:sum 2} is $0$. We use \eqref{E:inner sum} in \eqref{E:sum 2} and obtain
\begin{equation}
\label{E:sum 3}
\sum_{\substack{d|b\\(d,b')=1}}\mu(d)q^{r\left(\deg(b)-\deg(b')\right)-r\deg(d)}=q^{r\left(\deg(b)-\deg(b')\right)}\sum_{\substack{d|b\\(d,b')=1}}\mu(d)q^{-r\deg(d)}.
\end{equation}

We observe that the result we obtained in \eqref{E:sum 3} is independent of the particular choice of the tuple $(\alpha_1,\dots,\alpha_r)$. Now we compute, for the same fixed tuple $(\alpha_1,\dots,\alpha_r)$, the number of elements of $O$ of the form \eqref{E:element type} which are also in $U$. As explained before, using \eqref{E:big group}, we need to count the number of all elements of $U$ of the form \eqref{E:element type}, corresponding to this fixed tuple $(\alpha_1,\dots,\alpha_r)$, and satisfying \eqref{E:grad} and \eqref{E:gcd2}. 

For each $i\in\{1,\dots,r\}$, we let $a'_i(t):=t^{n_i}a_i(\frac{1}{t})$, where $a'_i(t)$ is a polynomial of degree less than $n_i$ (we are using the fact that $a_i(\frac{1}{t})\in\frac{1}{t}\mathbb{F}_q[\frac{1}{t}]$ is a polynomial of degree $n_i$). The requirement for an element of $O$ of the form \eqref{E:element type} to lie in $U$ is given by 
\begin{equation}
\label{E:U belong}
v_{\infty}\left(\frac{\alpha_i+b'q_i}{b}-\frac{a_i'}{t^{n_i}}\right)\ge n_i+1\text{ for each $i$.}
\end{equation}
Inequality \eqref{E:U belong} is equivalent with 
\begin{equation}
\label{E:U belong 2}
\deg\left(t^{n_i}(\alpha_i+b'q_i)-ba_i'\right)\le\deg(b)-1.
\end{equation}

Thus when we count the number of elements of $O\cap U$ which are of the form \eqref{E:element type} for one fixed tuple $(\alpha_1,\dots,\alpha_r)\in C$, we obtain the same sum as in \eqref{E:sum 1}, only that now we have the extra assumption \eqref{E:U belong 2} on top of the other restrictions on $q_i$. Again we can change the order of summation in \eqref{E:sum 1} and obtain \eqref{E:sum 2}, only that the inner summation is over all $q_i$ which also satisfy \eqref{E:U belong 2}. We obtain once again \eqref{E:q'} and we use it in \eqref{E:U belong 2} to get
\begin{equation}
\label{E:q' belong}
\deg\left(t^{n_i}(\alpha_i+b's_i+b'dq'_i)-ba_i'\right)\le\deg(b)-1.
\end{equation}

We are interested in estimating $\#(O\cap U)$ when $\deg(b)$ is much larger than the degree of $b'$ (and implicitly, much larger than the degrees of the $\alpha_i$, as $\deg(\alpha_i)<\deg(b')$) and also much larger than the numbers $n_i$. This is the case because our goal is to prove \eqref{E:simplified conclusion}, and as $k\rightarrow\infty$, $\deg(b)\rightarrow\infty$, while $U$ is fixed. Note that the degree of $b'$ is bounded by the degree of $P=\prod_{i=1}^lP_i^m$ (which is a fixed polynomial) and the numbers $n_i$ are fixed the moment we fixed $U$, while $b$ is the order of one of the elements of our infinite sequence of distinct torsion points. 

We go back now to \eqref{E:U belong 2}. We know that $\deg(q_i)<\deg(b)-\deg(b')$ (see \eqref{E:grad}). Thus, using \eqref{E:q'}, we get $\deg(q_i')<\deg(b)-\deg(b')-\deg(d)$. Therefore, let 
$$q'_i=\sum_{j=0}^{\deg(b)-\deg(b')-\deg(d)-1}\gamma^{(i)}_jt^j,$$
where some of the $\gamma^{(i)}_j\in\mathbb{F}_q$ could be $0$ (including some of the top coefficients). 

Let $L$ be the degree of $P$ and let $n_0:=\max\{n_1,\dots,n_r\}+1$. Assume for the moment that 
\begin{equation}
\label{E:degreed}
\deg(d)\le\deg(b)-L-n_0 <\deg(b)-\deg(b')-n_i,
\end{equation}
where in the second inequality from \eqref{E:degreed} we also used the fact that $b'|P$ (and so, $\deg(b')\le\deg(P)$). Under our assumption \eqref{E:degreed}, we obtain that the top $n_i$ coefficients $\gamma^{(i)}_j$ of $q'_i$ are determined by the coefficients of $d$, $b'$, $b$, $\alpha_i$ and $a_i'$ and by condition \eqref{E:q' belong}, while the remaining $(\deg(b)-\deg(b')-\deg(d)-n_i)$ coefficients $\gamma^{(i)}_j$ can be arbitrary elements of $\mathbb{F}_q$. Hence, under the assumption \eqref{E:degreed}, we obtain that the inner sum in \eqref{E:sum 2} associated only to the elements in $U$ of the form \eqref{E:element type} equals
\begin{equation}
\label{E:good degree}
q^{\sum_{i=1}^r\left(\deg(b)-\deg(b')-\deg(d)-n_i\right)},
\end{equation}
if $(b',d)=1$, while if $(b',d)>1$, the inner sum in \eqref{E:sum 2} is $0$.

Next we analyze the case $\deg(d)>\deg(b)-L-n_0$. In this case, \eqref{E:q' belong} shows that $\deg(q_i')<L+n_0$. Thus, the inner sum in \eqref{E:sum 2} can contribute at most $q^{r(L+n_0)}$ and this computation is without even taking into consideration the actual restrictions on the coefficients imposed by \eqref{E:q' belong}. 

Combining our findings in both cases with respect to assumption \eqref{E:degreed}, we conclude that the number of elements in $O\cap U$ which are of the form \eqref{E:element type} for a fixed $(\alpha_1,\dots,\alpha_r)\in C$ is
\begin{equation}
\label{E:sum 4}
\sum_{\substack{d|b\\ (d,b')=1\\ \deg(d)\le\deg(b)-L-n_0}}\mu(d)q^{r(deg(b)-\deg(b')-\deg(d))-\sum_{i=1}^r n_i}+\sum_{\substack{d|b\\ (d,b')=1\\ \deg(d)>\deg(b)-L-n_0}}O(q^{r(L+n_0)}),
\end{equation}
where the $O$-notation in the above second sum is the classical one and it refers in our context to the fact that the summand in the second sum from \eqref{E:sum 4} is at most $q^{r(L+n_0)}$ regardless of $b$. We note that the $O$-notation has nothing to do with our notation for orbits or for the set $O$, as the $O$-notation will always have attached to it, in parenthesis, a certain real number. 

Because of the degree condition for $d$ in the second sum in \eqref{E:sum 4}, we know there are at most $q^{L+n_0}$ possibilities for $d$ (for each fixed $b$), because $\deg\left(\frac{b}{d}\right)< L+n_0$. Hence, the second sum from \eqref{E:sum 4} is of the order of $q^{(r+1)(L+n_0)}$. Introducing this estimate in \eqref{E:sum 4} and adding and subtracting from \eqref{E:sum 4} the quantity
$$\sum_{\substack{d|b\\ (d,b')=1\\ \deg(d)>\deg(b)-L-n_0}}\mu(d)q^{r(deg(b)-\deg(b')-\deg(d))-\sum_{i=1}^r n_i},$$
we obtain the following estimate for \eqref{E:sum 4}:
\begin{equation}
\label{E:sum 5}
\sum_{\substack{d|b\\ (d,b')=1}}\mu(d)q^{r(deg(b)-\deg(b'd))-\sum_{i=1}^r n_i}-
 \sum_{\substack{d|b\\ (d,b')=1\\ \deg(d)>\deg(b)-L-n_0}} \mu(d)q^{r(deg(b)-\deg(b'd))-\sum_{i=1}^r n_i}+O\left(1\right).
\end{equation}
Note that in the $O$-estimate from \eqref{E:sum 5}, we replaced $O\left(q^{(r+1)(L+n_0)}\right)$ by $O(1)$, because $q$, $r$ and $L$ are always fixed, while $n_0$ is fixed the moment we fix $U$. 

We estimate the second sum in \eqref{E:sum 5} and we easily conclude it is also $O\left(q^{(r+1)(L+n_0)}\right)$. Hence the sum in \eqref{E:sum 5} is
\begin{equation}
\label{E:sum 6}
q^{r(\deg(b)-\deg(b'))-\sum_{i=1}^r n_i}\sum_{\substack{d|b\\ (d,b')=1}}\mu(d)q^{-r\deg(d)}+O\left(q^{(r+1)(L+n_0)}\right).
\end{equation}
The sum in \eqref{E:sum 6} is the number of elements in $O\cap U$ of the form \eqref{E:element type} for a fixed $(\alpha_1,\dots,\alpha_r)$. Note that $r,L,n_0$ are all constant as $\deg(b)\rightarrow\infty$. As explained before, the number of tuples $(\alpha_1,\dots,\alpha_r)\in C$ depends on $x_k$ (we recall that $O$ and $b$ were associated to some torsion element $x_k$ from our infinite sequence), but the cardinality of $C$ has the fixed upper bound $q^{Lr}$ because $b'|P$ and $\deg(P)=L$ and all of the polynomials $\alpha_i$ have degree less than $\deg(b')$ for each $i$. So, in order to prove \eqref{E:goal}, we use \eqref{E:sum 3} and \eqref{E:sum 6} and we are done if we show
\begin{equation}
\label{E:ultimate goal}
q^{r\left(\deg(b)-\deg(b')\right)}\sum_{\substack{d|b\\(d,b')=1}}\mu(d)q^{-r\deg(d)}\rightarrow\infty
\end{equation}
as $\deg(b)\rightarrow\infty$. Let $b_0$ be the product of all the powers of irreducible polynomials dividing $b$ other than the powers of the polynomials $P_1,\dots,P_l$. Thus $b_0$ is the largest divisor of $b$ coprime with $b'$. Then the sum in \eqref{E:ultimate goal} can be rewritten as
\begin{equation}
\label{E:sum 7}
\sum_{d|b_0}\mu(d)q^{-r\deg(d)}.
\end{equation}
The sum in \eqref{E:sum 7} equals
\begin{equation}
\label{E:sum 8}
\prod_{\substack{c\text{ irreducible}\\c|b_0}}\left(1-\frac{1}{q^{r\deg(c)}}\right).
\end{equation}
We observe that all the factors in the product \eqref{E:sum 8} are less than $1$ and so, if we extend the product \eqref{E:sum 8} to include also the possible prime divisors of $b$ from the set $\{P_1,\dots,P_l\}$, we can only decrease our product. So, to prove \eqref{E:ultimate goal}, it suffices to show
\begin{equation}
\label{E:final goal}
q^{r\left(\deg(b)-\deg(b')\right)}\prod_{\substack{c\text{ irreducible}\\c|b}}\left(1-\frac{1}{q^{r\deg(c)}}\right)\rightarrow\infty
\end{equation}
as $\deg(b)\rightarrow\infty$. We note that $\deg(b')\le L$ (as it was remarked previously in our proof, because $b'|P$ and $\deg(P)\le L$). So, we only need to show
\begin{equation}
\label{E:final goal 2}
q^{r\deg(b)}\prod_{\substack{c\text{ irreducible}\\c|b}}\left(1-\frac{1}{q^{r\deg(c)}}\right)\rightarrow\infty\text{ as }\deg(b)\rightarrow\infty.
\end{equation}
Consider the prime factorization of $b=\prod_{i=1}^sc_i^{e_i}$ in monic polynomials. Then the left hand side of \eqref{E:final goal 2} reads
\begin{equation}
\label{E:produs}
\prod_{i=1}^sq^{r(e_i-1)\deg(c_i)}\left(q^{r\deg(c_i)}-1\right).
\end{equation}
Because $\phi$ has modular transcendence degree at least $2$, the rank $r$ of $\phi$ is at least $2$ (otherwise, $\phi_t$ is a polynomial of degree $q$, which means that $\phi$ is isomorphic over $K^{\alg}$ with a Drinfeld module $\psi$ for which $\psi_t=t\tau^0+\tau$; hence $\phi$ would be isomorphic with a Drinfeld module defined over $F$). Because $r\ge 2$, we obtain $q^{r\deg(c_i)}-1\ge q^{\frac{r\deg(c_i)}{2}}$ for each $i$.

As $\deg(b)\rightarrow\infty$, then $\sum_{i=1}^se_i\deg(c_i)\rightarrow\infty$, which proves 
\begin{equation}
\label{E:final}
\prod_{i=1}^sq^{r(e_i-1)\deg(c_i)}\left(q^{r\deg(c_i)}-1\right)\ge q^{\frac{r\sum_{i=1}^se_i\deg(c_i)}{2}}\rightarrow\infty.
\end{equation}
This concludes the proof of \eqref{E:simplified conclusion}, which proves Theorem~\ref{T:equi}.
\end{proof}

\begin{remark}
Note that even if we used also in the last part of our argument the fact that $r\ge 2$, the only place where we used crucially the fact that $\phi$ has mdular transcendence degree at least $2$ was in Pink's result (Theorem~\ref{T:pink}). Therefore, if Pink's result were true for any Drinfeld module of generic characteristic, then our proof of Theorem~\ref{T:equi} would also hold for any such Drinfeld module. The limit in \eqref{E:final goal 2} is infinite even if $r=1$ (in which case $\phi$ would have modular transcendence degree $1$).
\end{remark}

\section{The Bogomolov Conjecture and the Manin-Mumford Theorem for Drinfeld modules}
\label{se:extensions}

In this section, all subvarieties are \emph{closed} subvarieties.

Let $\phi:A\rightarrow K\{\tau\}$ be a Drinfeld module of generic characteristic. For each positive integer $g$, we let $\phi$ act on $\mathbb{G}_a^g$ diagonally. Therefore, we may define just as before the torsion points of the action of $\phi$ on $\mathbb{G}_a^g$ as all the $g$-tuples $(x_1,\dots,x_g)$ for which there exists a nonzero $a\in A$ such that $\phi_a(x_i)=0$ for each $i$. We believe a similar equidistribution result as our Theorem~\ref{T:equi} holds for the torsion points of $\mathbb{G}_a^g$. Before stating our conjecture, we require the following definitions.

\begin{definition}
\label{D:alg phi module}
An algebraic $\phi$-submodule of $\mathbb{G}_a^g$ is a $K^{\sep}$-algebraic subvariety of $\mathbb{G}_a^g$ which is invariant under the action of $\phi$.
\end{definition}

\begin{definition}
\label{D:torsion subvariety}
A torsion subvariety of $\mathbb{G}_a^g$ is a translate of an irreducible algebraic $\phi$-submodule of $\mathbb{G}_a^g$ by a torsion point. 
\end{definition}

\begin{definition}
\label{D:strict}
A sequence of points $\{x_k\}\subset\mathbb{G}_a^g(K^{\alg})$ is strict if any proper torsion subvariety of $\mathbb{G}_a^g$ contains $x_k$ for only finitely many values of $k$.
\end{definition}

For a point $x\in\mathbb{G}_a^g(K^{\sep})$, we let as before $O(x)$ denote the (finite) orbit of $x$ under the diagonal action of $\aut$ on $\mathbb{G}_a^g(K^{\sep})$. We also define the associated probability measure $\dirac_x$ on $\C^g$ for such an orbit $O(x)$:
$$\dirac_x:=\frac{1}{\#O(x)}\sum_{y\in O(x)}\delta_y.$$
Finally, we denote by $\nu^{(g)}$ the product measure on $T^g$ corresponding to $\nu$ taken $g$ times.

\begin{conjecture}
\label{C:equi}
Let $\phi$ be a Drinfeld module of generic characteristic. Assume $\End(\phi)=A$. Let $\{x_k\}$ be a strict sequence of torsion points in $\mathbb{G}_a^g(K^{\alg})$, for some $g\ge 1$. Then $\dirac_{x_k}\goes\nu^{(g)}$.
\end{conjecture}

\begin{remarks}
We used in Conjecture~\ref{C:equi} the same convention as explained in Remark~\ref{R:support}, regarding the measures $\dirac_{x_k}$. Because their support is contained in $T^g$, we interpret them as probability measures on $T^g$, rather than as probability measures on the larger space $\C^g$.

We did not include in our Conjecture~\ref{C:equi} the hypothesis that $\phi$ has modular transcendence degree at least $2$ (as we did in our Theorem~\ref{T:equi}), because, as mentioned before, it is believed that Pink's Theorem~\ref{T:pink} holds without this extra hypothesis on $\phi$.

We require in our Conjecture~\ref{C:equi} the hypothesis on the sequence $\{x_k\}$ being strict because otherwise the support of the limit measure would lie on a proper subvariety of $T^g$. In the case $g=1$, our hypothesis in Theorem~\ref{T:equi} that the sequence $\{x_k\}\subset K^{\alg}$ contains distinct torsion points suffices for the condition that the sequence is strict (because the torsion subvarieties of $\mathbb{G}_a$ are the torsion points of $\phi$). Actually, our proof of Theorem~\ref{T:equi} follows precisely the same way under the slightly weaker hypothesis that the sequence $\{x_k\}\subset\phi_{\tor}$ contains each torsion point of $\phi$ at most finitely many times (this condition being equivalent with the condition that $\{x_k\}$ is strict). 
\end{remarks}

A positive answer to our Conjecture~\ref{C:equi} would provide a proof to the following result (the Manin-Mumford Theorem for Drinfeld modules of generic characteristic).

\begin{theorem}
\label{T:man-mum}
Let $\phi:A\rightarrow K\{\tau\}$ be a Drinfeld module of generic characteristic. Assume $\End(\phi)=A$. Let $g\ge 1$ and let $X$ be an irreducible $K^{\sep}$-subvariety of $\mathbb{G}_a^g$ (i.e., an irreducible closed subvariety of the $g$-dimensional affine space). If $X(K^{\sep})\cap\phi_{\tor}^g$ is Zariski dense in $X$, then $X$ is a torsion subvariety of $\mathbb{G}_a^g$.
\end{theorem}

As mentioned in Section~\ref{se:intro}, Theorem~\ref{T:man-mum} was proved by Scanlon in \cite{TS} using the methods of model theory of difference fields. His result is valid even without the extra assumption that the endomorphism ring of $\phi$ equals $A$. However, a positive answer to Conjecture~\ref{C:equi} would provide a completely different proof of Theorem~\ref{T:man-mum}, given purely in the language of number theory.

Moreover, we believe that an equidistribution result, similar to the results proved by Bilu, Ullmo and Zhang, holds also for points of small height associated to the action of a Drinfeld module (see \cite{T} or \cite{bogo} for the definition of the height $\hhat$ associated to a Drinfeld module). 

\begin{conjecture}
\label{C:small points}
Let $\phi:A\rightarrow K\{\tau\}$ be a Drinfeld module of generic characteristic and let $g\ge 1$. Assume $\End(\phi)=A$. Let $\{x_k\}\subset\mathbb{G}_a^g(K^{\sep})$ be a strict sequence. If $\lim_{k\rightarrow\infty}\hhat(x_k)=0$, then $\dirac_{x_k}\goes\nu^{(g)}$.
\end{conjecture}

\begin{remarks}
The measures $\dirac_{x_k}$ are probability measures on $\C^g$, while $\nu^{(g)}$ is the normalized Haar measure on $T^g$. Therefore, we interpret the conclusion of Conjecture~\ref{C:small points} as follows: the measures $\dirac_{x_k}$ converge weakly to the probability measure $\nu^{(g)}$ on $\C^g$, which is supported on $T^g$ (and the restriction of $\nu^{(g)}$ on $T^g$ is a Haar measure).

Conjecture~\ref{C:equi} is a particular case of Conjecture~\ref{C:small points} because all the torsion points of a Drinfeld module have height $0$.
\end{remarks}

An equidistribution result as Conjecture~\ref{C:small points} would lead to the following form of the Bogomolov Conjecture in the context of Drinfeld modules of generic characteristic.

\begin{conjecture}
\label{C:Bogo}
Let $\phi:A\rightarrow K\{\tau\}$ be a Drinfeld module of generic characteristic. Assume $\End(\phi)=A$. Let $g\ge 1$ and let $X$ be an irreducible $K^{\sep}$-subvariety of $\mathbb{G}_a^g$. For each $n\ge 1$, we let $$X_n:=\{x\in X(K^{\sep})\mid \hhat(x)<\frac{1}{n}\}.$$
If for each $n\ge 1$, $X_n$ is Zariski dense in $X$, then $X$ is a torsion subvariety of $\mathbb{G}_a^g$.
\end{conjecture}

We can also formulate the Manin-Mumford and the Bogomolov questions for Drinfeld modules of finite characteristic. However, we cannot always expect the conclusion be that the variety $X$ is a torsion subvariety (see Section $4$ of \cite{bogo} for a comprehensive discussion of the pathologies appearing in finite characteristic Drinfeld modules).

We will show how to deduce Theorem~\ref{T:man-mum} and Conjecture~\ref{C:Bogo} from Conjecture~\ref{C:equi} and, respectively Conjecture~\ref{C:small points}. More precisely, we will show that knowing the validity of Conjectures~\ref{C:equi} and \ref{C:small points} for all $g\le N$ implies the conclusions of Theorem~\ref{T:man-mum} and, respectively Conjecture~\ref{C:Bogo} for $g=N$. Our proofs are inspired by the arguments from the proof of Theorem $5.1$ of \cite{Bilu}.

We first need a relative notion of the condition that a sequence is strict.
\begin{definition}
\label{D:relative strict}
Let $Y$ be an algebraic $\phi$-submodule of $\mathbb{G}_a^N$ and let $\{x_k\}\in Y(K^{\sep})$ be a sequence of points in $Y$. We call the sequence $\{x_k\}$ strict relative to $Y$ if any torsion subvariety of $Y$ of dimension smaller than $\dim(Y)$ contains $x_k$ for only finitely many values of $k$.
\end{definition}

We will also use the following result.
\begin{lemma}
\label{L:components}
Let $Y$ be a subvariety of $\mathbb{G}_a^g$ such that $\phi_t(Y)\subset Y$. Then every irreducible component of $Y$ is a torsion subvariety.
\end{lemma}

\begin{proof}
Let $Y_1,\dots,Y_m$ be the irreducible components of $Y$. Then for each $i\in\{1,\dots,m\}$ and for each $n\ge 1$, there exists $j(i,n)\in\{1,\dots,m\}$ such that $\phi_{t^n}(Y_i)=Y_{j(i,n)}$ (because $\phi_{t^n}(Y_i)$ is also an irreducible component of $Y$). Hence, for each $i\in\{1,\dots,m\}$, there exist positive integers $n_1<n_2$ 9depending on $i$) such that $\phi_{t^{n_1}}(Y_i)=\phi_{t^{n_2}}(Y_i)$. Thus, $\phi_{t^{n_1}}(Y_i)$ is invariant under $\phi_{t^{n_2-n_1}}$. Using Lemme $4$ of \cite{Den} we conclude $\phi_{t^{n_1}}(Y_i)$ is a torsion subvariety. Therefore, using that $Y_i$ is irreducible, we obtain that also $Y_i$ is a torsion subvariety.
\end{proof}

The following lemma is a classical result, whose proof we include for completeness.
\begin{lemma}
\label{L:infinite sets are dense}
Let $S$ be an infinite set in $K^{\alg}$ and let $n\ge 1$. Then $S^n$ (the cartesian product of $S$ with itself $n$ times) is Zariski dense in $\mathbb{G}_a^n$.
\end{lemma}

\begin{proof}
We prove the statement of our lemma by induction on $n$. For $n=1$, the statement is clear, as every infinite set is Zariski dense in the $1$-dimensional affine space. Next we assume the lemma holds for $n$ and we will prove it for $n+1$.

Let $f\in K^{\alg}[X_1,\dots,X_{n+1}]$ be a polynomial vanishing on $S^{n+1}$. We will prove $f=0$, which will show that indeed, $S^{n+1}$ is Zariski dense in $\mathbb{G}_a^{n+1}$. For each $\alpha\in S$, $f(X_1,\dots,X_n,\alpha)$ vanishes on $S^n\subset\mathbb{G}_a^n$. Using the induction hypothesis, we conclude 
\begin{equation}
\label{E:vanishing}
f(X_1,\dots,X_n,\alpha)=0.
\end{equation}
We consider $f\in K^{\alg}(X_1,\dots,X_n)[X_{n+1}]$ as a polynomial of only the variable $X_{n+1}$. Because \eqref{E:vanishing} holds for the (infinitely many) elements $\alpha\in S$, we conclude $f=0$, as desired.
\end{proof}
Because $\phi_{\tor}\subset K^{\alg}$ is infinite, we obtain the following consequence of Lemma~\ref{L:infinite sets are dense}.
\begin{corollary}
\label{C:torsion is dense}
For each $n\ge 1$, the torsion submodule $\phi_{\tor}^n$ is Zariski dense in $\mathbb{G}_a^n$.
\end{corollary}
Moreover, the following is also true.
\begin{corollary}
\label{C:torsion is dense 2}
Let $Y\subset\mathbb{G}_a^N$ be an algebraic $\phi$-submodule. Then $\phi_{\tor}^N\cap Y$ is Zariski dense in $Y$.
\end{corollary}

\begin{proof}
Let $Y_0$ be the connected component of $Y$. Then $Y_0$ is isomorphic to $\mathbb{G}_a^M$ for $M:=\dim(Y)$. Using Corollary~\ref{C:torsion is dense}, we conclude $\phi_{\tor}^N\cap Y_0$ is Zariski dense in $Y_0$. Moreover, because all the irreducible components of $Y$ are translates of $Y_0$ by torsion points (see Lemma~\ref{L:components}), we conclude that $\phi_{\tor}^N\cap Y$ is Zariski dense in $Y$.
\end{proof}

We first show that the validity of Conjecture~\ref{C:equi} for all $g\le N$ yields the following key result.
\begin{theorem}
\label{T:max torsion subvarieties}
Let $Y$ be an algebraic $\phi$-submodule of $\mathbb{G}_a^N$ and let $X$ be a $K^{\sep}$-subvariety of $Y$. Then $X$ has at most finitely many maximal torsion subvarieties.
\end{theorem}

\begin{proof}[Proof of Theorem~\ref{T:max torsion subvarieties}.]
We prove our theorem by induction on $\dim(Y)$. The case $\dim(Y)=0$ is trivial, as then $Y$ consists of only finitely many torsion points.

Assume Theorem~\ref{T:max torsion subvarieties} holds for $\dim(Y)<M\le N$ and we will prove it for varieties of dimension $M$.

First we note that without loss of generality we may assume $\dim(X)<M$. Indeed, if $\dim(X)=M$, then the $M$-dimensional irreducible components of $X$ are also irreducible components of $Y$. But the irreducible components of $Y$ are torsion subvarieties by Lemma~\ref{L:components}. Removing the irreducible components of $X$ which are also irreducible components of $Y$ would make $X$ have strictly smaller dimension than $M$ and would not change the conclusion of Theorem~\ref{T:max torsion subvarieties}, because we would remove only finitely many maximal torsion subvarieties of $X$. Hence, we may assume $\dim(X)<M$.

Secondly, we may replace $X$ with $$\bigcup_{\sigma\in\aut}\sigma(X)$$ and replace $Y$ with $$\bigcup_{\sigma\in\aut}\sigma(Y),$$ which remains an algebraic $\phi$-module because the action of $\phi$ is invariant under $\aut$. If our Theorem~\ref{T:max torsion subvarieties} would fail for $X\subset Y$ in the first place, then it would also fail for the above varieties which replace them. The advantage of our reduction is that both $Y$ and $X$ are now invariant under $\aut$. Note that while making this reduction we do not change the dimension of $X$. So, we still have $\dim(X)<M$.

Assume $X$ has infinitely many maximal torsion subvarieties and let $Y_i$ be a complete list of them. For distinct $i$ and $j$, $Y_i\cap Y_j$ is a proper torsion subvariety for both $Y_i$ and $Y_j$ (as both $Y_i$ and $Y_j$ are maximal torsion subvarieties of $X$). Moreover, because torsion subvarieties are irreducible by definition, we conclude
\begin{equation}
\label{E:dimension}
\dim(Y_i\cap Y_j)<\min\{\dim(Y_i),\dim(Y_j)\}.
\end{equation}
Hence, using \eqref{E:dimension} we obtain that for each $k$, 
$$\left(\cup_{i=1}^{k-1}Y_i\right)\cap Y_k$$
is a proper torsion subvariety of $Y_k$. Because the torsion submodule $\phi_{\tor}^N$ is Zariski dense in $Y_k$ (see Corollary~\ref{C:torsion is dense 2}), we conclude that there exists a torsion point $x_k\in Y_k\setminus\left(\bigcup_{i=1}^{k-1}Y_i\right)$.

We will prove next that the sequence of torsion points $\{x_k\}$ is strict relative to $Y$. Let $Z$ be a torsion subvariety of $Y$ with $\dim(Z)<\dim(Y)$. We will show $X\cap Z$ contains finitely many $x_k$ (note that by construction, $\{x_k\}\subset X$).

Because $\dim(Z)<M$, we can apply the induction hypothesis and conclude there are finitely many maximal torsion subvarieties of $X\cap Z$. Indeed, let $Z=\alpha+W$, where $\alpha$ is a torsion point and $W$ is an irreducible algebraic $\phi$-submodule. Then 
\begin{equation}
\label{E:intersection}
X\cap Z=\alpha+\left(\left(-\alpha+X\right)\cap W\right).
\end{equation}
Thus we apply the inductive hypothesis to $X':=\left(-\alpha+X\right)\cap W$ and derive that $X'\subset W$ contains finitely many maximal torsion subvarieties. Therefore, $X\cap Z=\alpha+X'$ contains finitely many maximal torsion subvarieties. Let $W_1,\dots,W_l\subset X\cap Z$ be a complete list of them. Because they are torsion subvarieties contained in $X$, for each $i\in\{1,\dots,l\}$, there exists $j$ such that $W_i\subset Y_j$. But each $Y_j$ contains only finitely many $x_k$, by the construction of $\{x_k\}$. Hence, each of $W_i$ contains only finitely many of the $x_k$ and thus, $Z$ contains finitely many of $\{x_k\}$ (we recall that $\{x_k\}\subset X$, by construction). This proves that the sequence $\{x_k\}$ is strict relative to $Y$.

Because $\dim(Y)=M$, there exists a suitable projection $\pi$ of $Y$ on $M$ of the $N$ coordinates of $\mathbb{G}_a^N$ such that $\pi$ is a dominant morphism. At the expense of relabelling the coordinates of $\mathbb{G}_a^N$, we may assume $\pi:Y\rightarrow\mathbb{G}_a^M$. Moreover, because $Y$ is an algebraic group, $\pi(Y)$ is also an algebraic group. Using $\dim(\pi(Y))=M$ (because $\pi$ is a dominant morphism), we conclude $\pi(Y)=\mathbb{G}_a^M$. Because $\dim(Y)=M=\dim(\pi(Y))$, we conclude each fiber of $\pi$ is finite.

We claim the sequence $\{\pi(x_k)\}\subset\mathbb{G}_a^M(K^{\sep})$ is strict. Assume there exists some proper torsion subvariety $Z:=\alpha+W\subset\mathbb{G}_a^M$ which contains infinitely many $\pi(x_k)$ ($\alpha$ is a torsion point and $W$ is a proper algebraic $\phi$-submodule of $\mathbb{G}_a^M$). Let $S$ be the finite orbit of $\alpha$ under the action of $\phi$ on $\mathbb{G}_a^M$. Let $Z_0:=\cup_{\beta\in S}\left(\beta+W\right)$. Clearly, $Z_0$ is a proper algebraic $\phi$-submodule of $\mathbb{G}_a^M$ ($\dim(Z_0)=\dim(W)<M$). Moreover, by our assumption, $Z_0$ contains infinitely many $\pi(x_k)$.

Let $Z':=\pi^{-1}(Z_0)\subset Y$. Because $Y$ is invariant under $\phi_t$, then $\phi_t(Z')\subset Y$. Moreover, $$\pi(\phi_t(Z'))=\phi_t(\pi(Z'))=\phi_t(Z_0)=Z_0.$$
Hence $\phi_t(Z')\subset Z'$. Using Lemma~\ref{L:components} we conclude $Z'$ is a finite union of torsion subvarieties.

Because the kernel of $\pi$ is finite, $\dim(Z')=\dim(Z_0)<M=\dim(Y)$. Moreover, because $Z_0$ contains infinitely many $\pi(x_k)$, $Z'$ contains infinitely many $x_k$. Hence $Z'$ is a finite union of torsion subvarieties of $Y$ of dimension smaller than $M$, and $Z'$ contains infinitely many $x_k$. This contradicts our proof that $\{x_k\}$ is strict relative to $Y$. We conclude $\{\pi(x_k)\}$ is a strict sequence  of torsion points in $\mathbb{G}_a^M(K^{\sep})$.

Using Conjecture~\ref{C:equi} for $g=M\le N$ and for the strict sequence $\{\pi(x_k)\}\subset\phi_{\tor}^M$, we conclude $\dirac_{\pi(x_k)}\goes\nu^{(M)}$. By the second reduction step for our proof of Theorem~\ref{T:max torsion subvarieties}, $X$ is invariant under $\aut$. Thus $\pi(X)$ is invariant under $\aut$. Hence the measures $\dirac_{\pi(x_k)}$ are all supported on $\pi\left(X(K^{\sep})\right)$. But $\pi(X(K^{\sep}))\subset\pi(X(K^{\alg}))$, which is a closed set in the $v_{\infty}$-adic topology. Therefore, the weak limit $\nu^{(M)}$ is supported also on $\pi\left(X(K^{\alg})\right)$. But, by construction the support of $\nu^{(M)}$ is $T^M$, which contains the torsion submodule of $\mathbb{G}_a^M$. Therefore, $\pi(X)$ contains the torsion submodule of $\mathbb{G}_a^M$. As this torsion submodule is Zariski dense in $\mathbb{G}_a^M$ (see Corollary~\ref{C:torsion is dense}), we conclude $\pi(X)=\mathbb{G}_a^M$. Thus $\dim(X)=M=\dim(Y)$, which contradicts our first reduction step: $\dim(X)<M$. Therefore $X$ has finitely many maximal torsion subvarieties.
\end{proof}

Theorem~\ref{T:man-mum} for $g=N$ is an immediate corollary of Theorem~\ref{T:max torsion subvarieties}.
\begin{proof}[Proof of Theorem~\ref{T:man-mum}.]
Assume $X$ is not a torsion subvariety of $\mathbb{G}_a^N$. By Theorem~\ref{T:max torsion subvarieties} applied to $X\subset\mathbb{G}_a^N$, $X$ is not the union $Z$ of its maximal torsion subvarieties, because there are finitely many of them and each one has smaller dimension than $X$ (here we use the irreduciblity of $X$). By construction, $Z$ contains all the torsion points of $X$, which thus contradicts the hypothesis that the set of torsion points of $\mathbb{G}_a^N$ is dense in $X$. Therefore, $X$ is indeed a torsion subvariety of $\mathbb{G}_a^N$.
\end{proof}

Assuming the validity of Conjecture~\ref{C:small points} for all $g\le N$, we prove the following generalization of Conjecture~\ref{C:Bogo} for $g=N$.
\begin{theorem}
\label{T:star}
Let $Y$ be an algebraic $\phi$-submodule of $\mathbb{G}_a^N$. Let $X$ be a $K^{\sep}$-subvariety of $Y$ and let $Z$ be the union of all (finitely many) maximal torsion subvarieties of $X$. If $Z\ne X$, then there exists a positive constant $C$ (depending on $X$) such that for each $x\in \left(X\setminus Z\right)(K^{\sep})$, $\hhat(x)\ge C$.
\end{theorem}

We first note that because we assumed the validity of Conjecture~\ref{C:small points} for all $g\le N$, we also assume the validity of Conjecture~\ref{C:equi} for all $g\le N$, because Conjecture~\ref{C:equi} is a particular case of Conjecture~\ref{C:small points}. Hence Theorem~\ref{T:max torsion subvarieties} holds and we \emph{do} know that $X$ has finitely many maximal torsion subvarieties.

Before proving Theorem~\ref{T:star}, we sketch the proof of Conjecture~\ref{C:Bogo} using the result of Theorem~\ref{T:star} applied to $X\subset Y=\mathbb{G}_a^g$. If $X$ is not a torsion subvariety, then it is not equal with the finite union $Z$ of its maximal torsion subvarieties, because $\dim(Z)<\dim(X)$ (we also use here the fact that $X$ is irreducible in Conjecture~\ref{C:Bogo}). Hence, there exists $C>0$ as in Theorem~\ref{T:star}. Let $n$ be a positive integer such that $\frac{1}{n}<C$. Then $X_n$ (defined as in Conjecture~\ref{C:Bogo}) is a subset of $Z$, which contradicts the hypothesis that $X_n$ is dense in $X$. Therefore, our assumption was false and so, indeed $X$ is a torsion subvariety.

\begin{proof}[Proof of Theorem~\ref{T:star}.]
We proceed by induction on $\dim(Y)$. If $\dim(Y)=0$, then both $Y$ and $X$ are finite unions of torsion points. Hence the theorem is vacuously true. We assume Theorem~\ref{T:star} holds for $\dim(Y)<M\le N$ and we prove it also holds for $\dim(Y)=M$.

We may assume without loss of generality that $\dim(X)<M$. Otherwise, the irreducible components of $X$ of dimension $M$ are also irreducible components of $Y$ and so, by Lemma~\ref{L:components}, they are torsion subvarieties. Therefore, they are contained in $Z$. So, removing them will not change $X\setminus Z$.

At the expense of replacing $Y$ and $X$ by the respective finite unions of their orbits under the action of $\aut$, we may assume that both $Y$ and $X$ are invariant under $\aut$. Note that replacing $X$ with $\cup_{\sigma\in\aut}\sigma(X)$, replaces $X\setminus Z$ with $\cup_{\sigma\in\aut}\sigma\left(X\setminus Z\right)$, which contains $X\setminus Z$. 

Assume Theorem~\ref{T:star} does not hold for $X\subset Y$. Then we can find a sequence $\{x_k\}\subset \left(X\setminus Z\right)(K^{\sep})$ such that $\hhat(x_k)<\frac{1}{k}$ for each positive integer $k$. We prove next that the sequence $\{x_k\}$ is strict relative to $Y$.

Let $Y'$ be a torsion subvariety of $Y$ of dimension smaller than $M$ and assume $Y'$ contains infinitely many $x_k$. Let $Y'=\alpha+W$, where $\alpha$ is a torsion point and $W$ is an (irreducible) algebraic $\phi$-submodule of $Y$ of dimension smaller than $M$. Then, as in \eqref{E:intersection}, $$X\cap Y'=\alpha+\left(\left(-\alpha+X\right)\cap W\right)$$ and so, we can apply the inductive hypothesis to $X\cap Y'$ (because $\dim(W)=\dim(Y')<M$). Note that the height function is not changed under translations by torsion points (this allows us to pass the inductive hypothesis from $\left(-\alpha+X\right)\cap Y'$ to $X\cap Y'$).

We conclude that either $X\cap Y'$ equals the finite union $Z'$ of its maximal torsion subvarieties, or there exists a constant $C'>0$ such that for every $$x'\in\left(\left(X\cap Y'\right)\setminus Z'\right)(K^{\sep}),$$ $\hhat(x')\ge C'$. But the maximal torsion subvarieties of $X\cap Y'$ are contained in the maximal torsion subvarieties of $X$, which means that $Z'$ contains no points from the sequence $\{x_k\}\subset\left(X\setminus Z\right)(K^{\sep})$. On the other hand, if $X\cap Y'\ne Z'$, then $$\hhat(x_k)<\frac{1}{k}<C',$$ for every positive integer $k>\frac{1}{C'}$. So, there are finitely many points of the sequence $\{x_k\}$ contained in $\left(\left(X\cap Y'\right)\setminus Z'\right)(K^{\sep})$. Therefore, $Y'$ contains only finitely many points of the sequence $\{x_k\}\subset X$, which shows that indeed, $\{x_k\}$ is a strict sequence relative to $Y$. 

We can redo now the argument from the last part of the proof of Theorem~\ref{T:max torsion subvarieties}. Indeed, we find again a suitable projection $\pi:Y\rightarrow\mathbb{G}_a^M$ and prove as we did before, that $\{\pi(x_k)\}$ is a strict sequence (using that $\{x_k\}$ is a relative strict sequence). Because the height of a point in the affine space is the sum of the heights of each of its coordinates, $\hhat(\pi(x_k))\le\hhat(x_k)$, for each $k$. Hence $\lim_{k\rightarrow\infty}\hhat(\pi(x_k))=0$. Thus, we can apply the conclusion of Conjecture~\ref{C:small points} for the strict sequence $\{\pi(x_k)\}\subset\mathbb{G}_a^M(K^{\sep})$ of small points and conclude that $\dirac_{\pi(x_k)}\goes\nu^{(M)}$. This shows that the support $T^M$ of $\nu^{(M)}$ is contained in $\pi(X(K^{\alg}))$. Hence $\pi(X)=\mathbb{G}_a^M$ and so, $\dim(X)=\dim(Y)$, contradicting our assumption that $\dim(X)<M$. Therefore there exists a uniform positive lower bound $C$ for the height of the points in $\left(X\setminus Z\right)(K^{\sep})$.
\end{proof}

\end{document}